\documentclass[12pt]{article}

\usepackage[a4paper,left=1in,top=1in,width=6.3in,height=9in]{geometry}

\usepackage{amsmath,amsthm,amssymb,graphicx,bm}
\usepackage[sort]{natbib}
 \bibpunct{(}{)}{;}{a}{,}{,}

\newcommand{\tr}{\text{tr}}
\newcommand{\atan}{\text{atan}}
\newcommand{\nil}[1]{}
\newcommand{\la}{\langle}
\newcommand{\ra}{\rangle}
\newtheorem{thm}{Theorem}

\newcommand{\diag}{\text{diag}}
\newcommand{\sme}{\text{SME}}
\newcommand{\hsme}{\text{SME,hybrid}}
\newcommand{\mle}{\text{MLE}}

\newcommand{\var}{\text{var}}
\newcommand{\cov}{\text{cov}}
\newcommand{\acos}{\text{arccos}\; }
\begin{document}
\title{Score matching estimators for directional distributions}
\author{Kanti V Mardia, John T Kent, and Arnab K Laha \\ [5mm]
 University of Leeds and University of Oxford, University of Leeds,\\
 Indian Institute of Management Ahmedabad }
\date{}
\maketitle

\abstract{One of the major problems for maximum likelihood estimation
  in the well-established directional models is that the normalising
  constants can be difficult to evaluate. A new general method of
  ``score matching estimation'' is presented here on a compact
  oriented Riemannian manifold. Important applications include von
  Mises-Fisher, Bingham and joint models on the sphere and related
  spaces.  The estimator is consistent and asymptotically normally
  distributed under mild regularity conditions.  Further, it is easy to compute
  as a solution of a linear set of equations and requires no knowledge
  of the normalizing constant.  Several examples are given, both
  analytic and numerical, to demonstrate its good performance.}

Some key words: Exponential family, Fisher-Bingham distribution,
Riemannian manifold, sphere, torus, von Mises distribution.

\section{Introduction} 
\label{sec:intro}
A novel ``score matching estimator'' was proposed in \citet{Hyv05,Hyv07a}
as an alternative to the maximum likelihood estimator.  For exponential
family models, a key advantage of this new estimator is that it avoids
the need to work with awkward normalizing constants.  A minor limitation is
that it requires the underlying distributions to have sufficiently
smooth probability densities.

The score matching estimator was originally developed for densities on
Euclidean space.  The extension to Riemannian manifolds was sketched
in \citet{Dawid-Lauritzen05} and \citet{Parry.etal12}, but without a
detailed analysis.  Here we give a systematic investigation on a
compact oriented Riemannian manifold $M$.  The main focus is on
exponential family models.  Key applications include the
Fisher-Bingham and multivariate von Mises distributions, which lie on
spheres and related spaces.

\nil{
Section \ref{sec:Riemann} reviews the relevant background on
Riemannian manifolds.  The motivation and definition of the score
matching estimator are given in Sections
\ref{sec:sme-derivation}--\ref{sec:sme}.  Sections
\ref{sec:sphere}--\ref{sec:dda} illustrate the methodology for
directional distributions on spheres and products of spheres.  In the
special case of the von Mises distribution on the circle it is
possible to give a more comprehensive investigation (Section
\ref{sec:efficiency}).  Some numerical examples are given in
Section \ref{sec:examples} to compare the score matching estimators to
other commonly used estimators.  Finally, Section \ref{sec:future}
gives a discussion.
}

\section{Background on Riemannian manifolds}
\label{sec:Riemann} 

A $p$-dimensional Riemannian manifold $M$, say, is characterized by a
metric tensor $G = G(x) = (g_{ij}(x))$, where $G(x)$ is a $p \times p$
symmetric positive definite matrix.  Here $x$ is a $p$-dimensional
vector representing a typical element of $M$ in local coordinates.
Let $G^{-1}(x) = (g^{ij}(x))$ denote the inverse matrix.  A uniform
measure on $M$ can be defined in local coordinates by $\mu(dx)=
\{\det(G(x))\}^{1/2} dx$.  See, e.g.,
\citet{Rosenberg97,Stein-Weiss71,Jost05} for a background on
Riemannian manifolds.

Next, let $u=u(x), \ v=v(x), \ x \in M$, be two real-valued functions
on $M$.  Let $\nabla u = (\partial u / \partial x_j, \ j=1, \ldots,
p)^T$ denote the gradient in local coordinates, i.e. the column vector
of partial derivatives.  An inner product on gradient vectors can be
defined by
\begin{equation}
\label{eq:ip}
\la u,v \ra = \la u,v \ra  (x)  = (\nabla u)^T G^{-1} (\nabla v),
\end{equation}
treated as a function of the local coordinate vector $x$.  

The Laplace-Beltrami operator, acting on $u$, is  defined
in local coordinates by
\begin{equation}
\label{eq:LB}
\Delta_M u = \sum_{i,j=1}^p \{\det(G(x))\}^{-1/2}\partial/\partial x_i 
\left[\{\det(G(x))\}^{1/2} g^{ij}(x) \partial u / \partial x_j \right].
\end{equation}
treated as a function of the local coordinate vector $x$.  Although
the gradient vector $\nabla u$ depends on the choice of local coordinates,
the uniform measure $\mu(dx)$, the gradient inner product
(\ref{eq:ip}) and the Laplace-Beltrami operator (\ref{eq:LB}) are
invariant under a change of local coordinates.

Stokes' Theorem, also known as the divergence theorem, connects the
gradient inner product and the Laplace-Beltrami operator.  Further, if
$M$ is assumed compact and oriented, there is no need for boundary
conditions.
\begin{thm}[Stokes' Theorem]
\label{th:stokes}
If $M$ is a compact oriented Riemannian manifold and $u(x),\ v(x)$ are
twice continuously differentiable functions on $M$, then
\begin{equation}
\label{eq:div}
\int_M \la u,v\ra (x) \, \mu(dx) = - \int_M (\Delta_M u) v\, \mu(dx).
\end{equation}
\end{thm}

If $M$ is isometrically embedded as a $p$-dimensional surface in a
Euclidean space, $M \subset \mathbb{R}^q$, then a point in $M$ can be
represented either as a $p$-dimensional vector $x$ in local
coordinates or as a $q$-dimensional vector $z = (z_1, \ldots, z_q)^T$,
say, in Euclidean coordinates. In this setting it is possible to give
simpler representations of the gradient inner product and the
Laplace-Beltrami operators.

Given a point $z$ in $M$, the tangent plane to $M$ at $z$ is a
$p$-dimensional hyperplane.  Let $P(z) = P = (p_{ij})$ denote the $q
\times q$ orthogonal projection matrix onto this hyperplane.  Then
$P=P^T, \ P=P^2$ and $P$ has rank $p$.  Suppose the functions $u$ and
$v$ have been extended to a neighbourhood $N$ of $M$ as functions
$\tilde{u}(z), \ \tilde{v}(z)$ of $z \in N \subset \mathbb{R}^q$, and
define the usual Euclidean gradient
\begin{equation*}
\label{eq:gradient-Eucl}
\nabla_E \tilde{u} = (\partial \tilde{u} / \partial z_j, \ j=1, \ldots,q)^T.
\end{equation*}

The gradient inner product and the Laplace-Beltrami operator were
given by (\ref{eq:ip}) and (\ref{eq:LB}) in local coordinates.  They can
also be written in Euclidean coordinates as
\begin{equation}
\label{eq:ip-Eucl}
\la u,v\ra = \left(\nabla_E \tilde{u}\right)^T P 
\left( \nabla_E \tilde{v} \right)
           = \left( P \nabla_E \tilde{u}\right)^T 
\left(P \nabla_E \tilde{v} \right)
\end{equation}
and
\begin{equation*}
\label{eq:LB-Eucl}
\Delta_M u = \tr \left\{ P \nabla_E^T 
\left( P \nabla_E  \tilde{u} \right) \right\}
= \sum_{i,j,k=1}^q  p_{ij} \partial/\partial z_i 
\left( p_{jk} \partial \tilde{u} / \partial z_k \right).
\end{equation*}
Note the derivatives in (\ref{eq:ip-Eucl}) are $q$-dimensional whereas
the derivatives in (\ref{eq:ip}) are $p$-dimensional; however, the
inner product is the same.

The simplest example is the sphere $S_2 = \{z \in \mathbb{R}^3 : z^Tz
= 1\}$.  In polar coordinates, $x = (\theta, \phi)^T$ with $\theta = $
colatitude and $\phi = $ longitude, the metric tensor becomes
$$
G = \begin{bmatrix} 1 & 0 \\ 0 & \sin^2 \theta \end{bmatrix},
$$
and the uniform measure becomes $\sin \theta d \theta d \phi$.
The Laplace-Beltrami operator becomes
$$
\Delta_M u = \partial^2 u / \partial \theta^2 + \cot \theta \partial u /
\partial \theta + (\sin \theta)^{-2} \partial^2 u / \partial \phi^2.
$$
The Euclidean embedding takes the form $z_1 = \cos \theta, \ z_2 =
\sin \theta \cos \phi, \ z_3 = \sin \theta \sin \phi$, and the
projection matrix is $P = I_3 - z z^T$.  It is straightforward to
check that gradient inner products in (\ref{eq:ip}) and
(\ref{eq:ip-Eucl}) are two ways of writing the same function.

\section{The score matching criterion 
on a compact oriented Riemannian manifold}
\label{sec:sme-derivation}
One way to motivate a statistical estimator is through the
minimization of a divergence between two probability distributions.
The most common example is the Kullback-Liebler divergence, which
leads to the maximum likelihood estimator.  In this paper we use a
divergence due to Hyv\"arinen, which leads to the score matching estimator.

Let $f$ and $f^*$ be two probability densities on a Riemannian
manifold $M$, defined with respect to the uniform measure $\mu$, where
$f$ and $f^*$ are assumed to be everywhere nonzero and twice
continuously differentiable.  Define the \emph{Hyv\"arinen divergence}
\citep{Hyv05,Hyv07a}
between the two densities in terms of an integrated gradient inner
product for the log ratio,
\begin{align} 
\Phi(f; f^*) &= \frac12 \int_M  \la \log (f/f^*), 
  \log (f/f^*) \ra\, f^*(x) \mu(dx) \notag\\
&= 
 \frac12 \int_M \left\{  \la \log f, \log f\ra 
-2 \la \log f, \log f^* \ra +   \la \log f^*, \log f^*\ra \right\}  
\, f^*(x)\, \mu(dx). 
\label{eq:H-div}
\end{align}
Hyv\"arinen proposed this divergence in the setting where $M$ is a Euclidean
space.  In this setting the metric tensor $G=I$ is hardly needed, and
the the gradient inner product (\ref{eq:ip}) simplifies to $\la u,v \ra =
\sum_{j=1}^p (\partial u / \partial x_j) (\partial v / \partial x_j)$.

For the rest of the paper we limit attention to the setting where $M$
is compact and oriented.  In some ways this setting is more
complicated than the Euclidean case.  The gradient inner product $\la
\cdot \ra$ is now given by (\ref{eq:ip}); it can also be written as
(\ref{eq:ip-Eucl}) when the manifold $M$ can be embedded in a
Euclidean space.  However, in other ways this setting is simpler.
Stokes' Theorem (\ref{eq:div}) holds automatically without the need for boundary
conditions and all continuous functions on $M$ are automatically
bounded and integrable.

Two simple properties ensure that minimizing (\ref{eq:H-div}) over $f$ is 
an identifiable approach to finding $f^*$.

\begin{description}
\item[Property 1.] If $f = f^*$, then $\Phi(f; f^*) = 0$.

\item[Property 2.] If $f \neq f^*$, then $\Phi(f; f^*) > 0$.
\end{description}

To prove these properties note by inspection in (\ref{eq:H-div}) that
$\Phi(f; f^*) \geq 0$ for all choices of densities and that $\Phi(f^*;
f^*) = 0$.  Conversely, if $f \neq f^*$, then we claim there must be
at least one point $x$ in some set of local coordinates such that
$\nabla f(x) \neq \nabla f^*(x)$; hence the integral (\ref{eq:H-div})
must be strictly positive.  For if there were no such $x$, a
contradiction would arise: the gradients would be everywhere equal and
so the densities would be the same up to a proportionality constant;
further, since both densities integrate to 1, the proportionality
constant would have to equal 1.

For fitting purposes, let $f^*(x)$ be regarded as the ``true''
distribution of the data and let $f(x) = f(x; \pi)$ denote a
parametric model, where $\pi$ is an $m$-dimensional vector of
parameters.  Then a ``best-fitting'' model can be defined by minimizing
$\Phi(f; f^*)$ over the parameters $\pi$.
 
Since the final term in (\ref{eq:H-div}) does not depend on $f$, it
can be dropped from the minimization criterion.  Further, by Stokes'
Theorem (\ref{eq:div}) and the simple result, $\partial \log f^*(x) / \partial x_j =
\{ f^*(x)\}^{-1} \partial f^*(x) / \partial x_j$, the middle term
becomes
\begin{align}
-\int_M \la \log f, \log f^* \ra  f^*(x)\, \mu(dx) 
&= -\int_M \la \log f, f^* \ra (f^*)^{-1} f^*(x)\, \mu(dx) \notag\\
&= -\int_M \la \log f, f^* \ra \,\mu(dx) \notag\\
&= \int_M (\Delta_M \log f) f^*(x)\, \mu(dx) \notag .
\end{align}

The final term in (\ref{eq:H-div}) does not depend on $f$.  Hence
minimizing $\Phi(f; f^*)$ over the parameters $\pi$ in $f$ is
equivalent to minimizing
\begin{align}
\Psi(f;f^*) &= \frac12 \int_M 
\left\{\la\log f, \log f\ra + 2 (\Delta_M \log f) \right\}
f^*(x) \mu(dx) \notag \\
&= \frac12E \left\{\la\log f, \log f\ra + 2 (\Delta_M \log f) \right\}, 
\label{eq:simplified-Fisher-divergence}
\end{align}
where $E$ denotes expectation under the measure $f^*(x) \mu(dx)$.
At this stage it is no longer necessary to impose any
regularity conditions on $f^*$; indeed it can be replaced by any
probability measure $F^*$, say, in which case we write $\Psi(f,F^*)$.

For the rest of the paper we specialize to the case where $f(x; \pi)$
forms a canonical exponential family on a compact oriented Riemannian
manifold $M$ with density
\begin{equation}
\label{eq:exp-fam}
f(x) = f(x; \pi)  \propto \exp\{\pi^T t(x)\}
\end{equation}
with respect to the uniform measure $\mu(dx)$, where $\pi$ is an
$m$-vector of natural parameters and $t(x) = (t_1(x), \ldots, t_m(x)$
is a vector of sufficient statistics.  The sufficient statistics are
assumed to satisfy the following regularity conditions.
\begin{itemize}
\item[A1] The constant function 1 and the functions $t_\ell(x) \
  (\ell=1, \ldots,m)$ are linearly independent on $M$.
\item[A2] The functions $t_\ell(x) \ (\ell=1, \ldots,m)$ are twice continuously
differentiable with respect to $x$.
\end{itemize}
The first assumption ensures the identifiability of $\pi$; different
values of $\pi$ correspond to different distributions.  The second
assumption justifies the use of Stokes' Theorem.

In this exponential family setting (\ref{eq:simplified-Fisher-divergence})
simplifies to
\begin{equation}
\label{eq:qf}
\Psi(f;F^*) = \frac12 \pi^T W \pi - \pi^T d,
\end{equation}
where $W = (w_{\ell_1\ell_2})$ and $d = (d_\ell)$ have elements
\begin{equation}
\label{eq:W-d-def}
w_{\ell_1\ell_2} = E \{ \la t_{\ell_1}, t_{\ell_2}\ra (x)\}, \quad
d_\ell = -E\{ \Delta_M t_\ell(x)\} \quad (\ell, \ell_1, \ell_2 = 1, \ldots, m).
\end{equation}

Minimizing (\ref{eq:qf}) yields the moment equation $W \pi - d = 0$.
In particular, if the ``true'' distribution $F^*(dx) = f(x; \pi_0)$ lies in
the parametric family, it follows that $\pi_0$ can be recovered from
the moments in $W$ and $d$ by
\begin{equation*}
\label{eq:moment-eq}
\pi_0 = W^{-1} d.
\end{equation*}
The identifiability property, Property 2  mentioned below
(\ref{eq:H-div}), implies that $W$ must be nonsingular.  In later
sections we illustrate how to calculate $W$ and $d$ in particular
cases.

In many examples the components of $t(x)$ are eigenfunctions of the
Laplace-Beltrami operator, so that $\Delta_M t_\ell(x) = -\lambda_\ell
t_\ell(x)$, for suitable constants $\lambda_\ell>0 \ (\ell=1, \ldots,
m)$ See, e.g., \citet[Ch 3.3]{Pat-Ell16} for the eigenfunctions of
various manifolds arising in Statistics; the specific example of the
sphere is discussed below.

\section{The score matching estimator and its properties}
\label{sec:sme}
Let $F^* = F_n$ denote the empirical distribution for a set of  data
$\{x_h, \ h=1, \ldots,n\}$.  The value $\hat{\pi}$ minimizing
$\Psi(f; F_n)$ is called the \emph{score matching estimator}.  More
explicitly,  
\begin{equation*}
\label{eq:sme}
\hat{\pi}_\sme = W_n^{-1} d_n,
\end{equation*}
where $W_n$ and $d_n$ are obtained from (\ref{eq:W-d-def}) after
replacing the expectation under $F^*$ by a sample average over the $n$
data points.  

To discuss the asymptotic sampling properties of the score matching
estimator, suppose the data comprise a random sample from the
exponential family model $f(x; \pi)$ with $\pi=\pi_0$.  Since $W$ is
nonsingular in the population case, it follows from the strong law of
large numbers, that $W_n$ must be positive definite with probability 1
for sufficiently large $n$.  This point has not generally been
emphasized in the literature.

Further, in many models this nonsingularity statement about  $W_n$ can be
strengthened to conclude that there is a fixed value $n_0$, depending on
the model but not on the data, such that $W_n$ must be positive
definite with probability 1 for $n \geq n_0$.  An analogous result for
the $p$-dimensional multivariate normal distribution states that the
sample covariance matrix is guaranteed to be positive definite with
probability 1 for $n \geq p+1$.

The limiting behaviour of $\hat{\pi}$ is straightforward to describe.
By the central limit theorem, $d_n$ and $W_n$ are asymptotically
jointly normally distributed with population means $d$ and $W$.  Hence
by the delta method, it follows that
$$
n^{1/2} \{\hat{\pi}_\sme - \pi_0\} \sim N_m(0, \Omega)
$$
for some limiting covariance matrix $\Omega$.  Further, by the
asymptotic optimality of maximum likelihood, $\Omega \geq
\mathcal{I}^{-1}$ under the usual ordering for positive semi-definite
matrices, where $\mathcal{I}$ denotes the Fisher information matrix.

The discussion here is limited to the exponential family case.
\citet{Forbes-Lauritzen14} note that extra regularity conditions are
needed for consistency and asymptotic normality when looking at score
matching estimation for more general densities.

The term ``score'' has several distinct
connotations in estimation.  (a) Conventionally, the ``score'' refers
to a derivative of the log likelihood with respect to the parameters;
it has close connections to maximum likelihood estimation.  (b)
However, in the context of the score matching estimator, the ``score''
refers to the derivative of the log likelihood with respect to the
state variable $x$.  (c) In addition, the term ``scoring rule''
\citep[e.g., ][]{Parry.etal12,Forbes-Lauritzen14} refers to a more
general function of $x$ and a distribution. Each scoring rule
determines a divergence, and the minimization of the divergence leads
to an estimator.  A scoring rule is different from the scores in (a)
and (b), though suitable choices for scoring rules lead to both the
maximum likelihood and score matching estimators.

\section{Details for the sphere}
\label{sec:sphere}
For this paper the most important choice for the manifold $M$ is the
unit sphere $S_p = \{z \in \mathbb{R}^q : \sum z_j^2 = 1\}$, a
$p$-dimensional manifold embedded in $\mathbb{R}^q, \ q=p+1$.  There
are two natural coordinate systems: embedded Euclidean coordinates $z$
and local coordinates $x$, i.e. polar coordinates in this case.

First we set out the key steps for the derivative calculations.  Let
$\tilde{u}(z)$ denote a scalar-valued function in Euclidean
coordinates.  The projected Euclidean gradient vector becomes
$$
P \nabla_E \tilde{u}  = (I-z z^T) \nabla_E  \tilde{u}, \quad P=I_q - z z^T.
$$

The eigenfunctions of $\Delta_M$ on $S_p, \ p \geq 1$, are known as
the spherical harmonics.  A spherical harmonic of degree $k \geq 0$
has eigenvalue $-\lambda_k$ where $\lambda_k = k(k+p-1)= k(k+q-2)$; see, e.g.,
\citet[p. 35]{Chavel84}, \citet[p. 125]{Pat-Ell16}.  The action of $\Delta_M$ 
on the linear
and quadratic spherical harmonics, expressed in Euclidean coordinates,
can be summarized as follows,
\begin{align}
\label{eq:LB-eigenfns}
\begin{split}
&\Delta_M z_j = -\lambda_1 z_j,\\
&\Delta_M (z_i^2 - z_j^2) = -\lambda_2 (z_i^2 - z_j^2),\\
&\Delta_M (z_i z_j) = -\lambda_2 (z_i z_j),
\end{split} \end{align} 
where $i \neq j \in \{1,\ldots, q\}$.  A systematic construction of
higher order spherical harmonics is given in
\citet[pp. 137--152]{Stein-Weiss71}, but they will not be needed
here. Statistical models involving spherical harmonics of degree
greater than two are straightforward in principle
\citep[e.g.][]{Beran79}, but in practice the components of $t(z)$ will
usually consist of linear and quadratic functions of $z$.
 
To illustrate these calculations, consider the Fisher-Bingham density
on $S_p$, 
\begin{equation}
\label{eq:fb}
f(z) \propto \exp \left\{b^Tz + z^T A z \right\},
\end{equation}
where $b$ is a $q$-vector, and $A$ is a $q \times q$ symmetric
matrix, $q=p+1$.  To ensure identifiability, the side condition
$\sum_{j=1}^{q} a_{jj} = 0$ is imposed.

The density can be recast as
\begin{equation*} 
\label{eq:fb-full}
\begin{split}
f(z) &\propto \exp \{ b_1 z_1 + \cdots + b_q z_q + a_{11} (z_1^2 - z_q^2) +
\cdots + a_{q-1,q-1} (z_{q-1}^2 - z_{q}^2) +\\
& \ \ \ \ \ \ a_{12} (2z_1z_2) + \cdots + a_{1q} (2z_1 z_{q}) + \cdots + 
a_{q-1,q} (2z_{q-1} z_{q}) \}\\
& = \exp \left\{ \sum_{\ell=1}^m \pi_\ell t_\ell(z) \right\}.
\end{split} \end{equation*}

Here the vector $\pi$ denotes the $m = q + q-1 + q(q-1)/2$
parameters in $b$ and $A$ in the order listed, and the vector
$t = t(z)$ denotes the corresponding functions of $z$; that is, 
the linear terms, the diagonal quadratic terms, and the cross-product
quadratic terms, respectively.

Next we gather the information needed to evaluate $W$ and $d$ in
(\ref{eq:W-d-def}).  For each sufficient statistic $t_\ell =
t_\ell(z)$, create a vector-valued function $u_\ell = u_\ell(z)
= \nabla_E t_\ell(z)$ by taking its Euclidean gradient, and create a 
scalar-valued function $v_\ell =
v_\ell(z) = z^T u_\ell$.  Then create an $m \times m$ matrix $W$
with entries
\begin{equation}
\label{eq:W-sphere}
w_{\ell_1 \ell_2} = E\{u_{\ell_1}^T u_{\ell_2} - v_{\ell_1} v_{\ell_2}\}.
\end{equation}
Table \ref{table:sphere}  gives the entries for  $u_\ell$
and $v_\ell$.  Also create a vector $d$ with entries
\begin{equation}
\label{eq:d-sphere}
d_\ell = -E\{ \Delta_M  t_{\ell}(z) \},
\end{equation}
using equation (\ref{eq:LB-eigenfns}) to evaluate the Laplace-Beltrami
operator.

\begin{table}
  \caption{Gradient details for the Fisher-Bingham density (\ref{eq:fb}) on
$S_q$.  Here $t_\ell(z)$ and $v_\ell(z)$ are scalars; $u_\ell(z)$ is a
$q$-dimensional vector. Also $e_j$ represents a unit vector along the
$j$th coordinate axis $(j=1, \ldots, q)$.}
\label{table:sphere}
\begin{center} \begin{tabular}{ccc}
$t_\ell(z)$ & $u_\ell(z)$ & $v_\ell(z)$ \\
$z_j$       & $e_j$       & $z_j$\\
$z_j^2 - z_p^2$ & $2(z_j e_j - z_q e_q)$ & $2(z_j^2 - z_q^2)$\\
$2z_i z_j$       & $2(z_i e_j + z_j e_i)$    & $4 z_i z_j$.
\end{tabular} \end{center}
\end{table}

\section{Hybrid estimators and reduced models}
\label{sec:hybrid}
The parameters for directional distributions can typically be split
into two parts: the concentration parameters and the orientation
parameters, where the normalizing constant depends just on the
concentration parameters.  Further it is often possible to estimate
the orientation parameters explicitly using sample moments.  Often
this orientation estimator can be viewed as an exact or an approximate
maximum likelihood estimator, and it can be computed without needing
estimates of the concentration parameters.  

Further, if the orientation parameters are known, then the distribution 
of the data becomes a natural exponential family for the concentration 
parameters.  We call this latter model \emph{reduced} because the number of
concentration parameters is smaller than the number of original parameters.

Hence, the following hybrid strategy provides a tractable estimation
procedure:

\begin{itemize}
\item[(a)] Split the parameters for the full model into orientation
and concentration parameters.

\item[(b)] Estimate the orientation parameters for the original
  data and the transform the data to a standardized form.

\item[(c)] After standardization, the estimated orientation parameters
  take a simple canonical form.  Treating them as known for the
  standardized data, the concentration parameters become natural
  parameters in a reduced exponential family.

\item[(d)] Use the score matching estimator to estimate the
  concentration parameters of the reduced model for the standardized
  data.
\end{itemize}

The resulting estimator can be called the \emph{hybrid score matching
  estimator} and denoted by $\hat{\pi}_\hsme$.  

The next section gives several examples from directional statistics to
illustrate this estimation strategy.  Table \ref{table:hybrid}
summarizes the form of each density.  See, e.g. \citet{Mardia-Jupp00,
  Jamm-Sen01} for background information on these distributions.  The
first three examples lie on the sphere $M=S_p$; the last example lies
on the torus, a direct product of $k$ circles, $M= (S_1)^k$.  In each
case the details (a)--(d)  are specified explicitly.

In all cases we assume a sample of size $n$ from the stated
distribution.  On the sphere the data are represented by an $n \times
q$ matrix $Z$ whose rows $z_h^T \ (h=1,\ldots, n)$, say, are
q-dimensional unit vectors in Euclidean coordinates.  The models are
special cases of the Fisher-Bingham density (\ref{eq:fb}) and Table
\ref{table:sphere} gives the details needed for the gradient
calculations.  On the torus, it is more convenient to use polar
coordinates to represent the data as an $n \times k$ matrix $\Theta$
of angles lying in $[0, 2\pi)$, with gradient calculations carried out
directly.  In each case the key step is to derive the formulas for the
matrix $W_n$ and the vector $d_n$ in the reduced model.

\begin{table}[t]
  \caption{Log densities for various standard directional distributions in 
    full and reduced form.  In full form, the models include both 
    orientation and concentration parameters.  In reduced form there 
    are just concentration parameters.  The models named in square 
    brackets are too general to have a useful reduced form, but the 
    submodels listed below them do have a useful reduced form.  The 
    first four models lie on the sphere $S_p$.  The final two models 
    are multivariate (MV) von Mises models lying on the torus $(S_1)^k$.}
\label{table:hybrid}
\renewcommand{\arraystretch}{1.2}
\begin{tabular}{ccc}
Name & Full & Reduced \\
von Mises-Fisher & $\mu^T z$ & $\kappa z_1$ \\
Bingham & $z^T A z$ & $z^T \Lambda z = \sum_{j=1}^{q-1} \lambda_j (z_j^2 - z_p^2)$ \\
$[$Fisher-Bingham$]$ & $\mu^T z + z^T A z$ & ---\\
Kent & 
$A \mu = 0$; 
$\lambda_1 = 0, \ \lambda_2 = -\lambda_3 = \beta$
& $\kappa z_1 + \beta(z_2^2 - z_3^2)$\\
$[$ MV von Mises$]$ & 
\parbox[c]{4cm}
{$$\sum_{r=1}^k \mu^{(r)T} z^{(r)} + \sum_{r<s} z^{(r)T} \Omega^{(r,s)} z^{(s)}$$} & 
---\\
MV von Mises sine & 
\parbox[c]{4cm}
{$\begin{cases} 
\sum_{r=1}^k \kappa^{(r)} c_r' + \sum_{r<s} \lambda^{(rs)} s_r' s_s',\\
c_r' = \cos \theta^{(r)\prime}, \ 
 s_r' = \sin \theta^{(r)\prime},\\
\theta^{(r)\prime} = \theta^{(r)} - \theta^{(r)}_0
\end{cases}$} & 
Set $\theta^{(r)}_0 = 0$.
\end{tabular}
\end{table}

\nil{ replaced entry in table
$\sum_{r=1}^k \kappa^{(r)} z_1^{(r)} + 
            \sum_{r<s} \lambda^{(rs)} z_2^{(r)} z_2^{(s)}$
}

\section{Directional  distributions}
\label{sec:dda}
\subsection{von Mises-Fisher distribution}
\label{sec:dda:vmf}
\begin{itemize}
\item[(a)] The density for the full von Mises-Fisher distribution
  takes the form in (\ref{eq:exp-fam}) with $t(z) = (z_1, \ldots,
  z_q)^T$ being the vector of linear functions in $z$.  This model
  forms a canonical exponential family.  The parameter vector can be
  written $\pi = \kappa \mu_0$ where $\kappa \geq 0$ is a scalar
  concentration parameter and $\mu_0$ is a unit orientation vector. On
  the circle, it is sometimes convenient to write $\mu_0 = (\cos
  \theta_0, \ \sin \theta_0)^T$ in polar coordinates.

\item[(b)] For the data matrix $Z (n \times q)$, the sufficient statistic is the
  $q$-dimensional sample mean vector $\overline{z}$.  The maximum
  likelihood estimate of $\mu_0$ is the unit vector
  $\hat{\mu}_{0,\mle} = \overline{z}/ ||\overline{z}||$ and is also
  the hybrid score matching estimator.  Let $R$ be a $q
  \times q$ orthogonal matrix such that $R^T \hat{\mu}_{0,\mle} =
  e_1$, where $e_1$ is a unit vector along the first coordinate axis
  in $\mathbb{R}^q$, and let $Y= Z R$, i.e.  $y_h = R^T z_h \ (h=1,
  \ldots,n)$, denote the standardized data.

\item[(c)] As shown in Table \ref{table:hybrid}, the reduced model for $Y$ 
involves just a single concentration parameter $\kappa$.

\item[(d)] For the reduced model, $W_n$ and $d_n$ in
  (\ref{eq:W-sphere})--(\ref{eq:d-sphere}) are one-dimensional, 
$$
W_n = 1-\frac1n \sum y_{h1}^2, \quad d_n = (q-1) \sum y_{h1},
$$
so that the hybrid score matching estimator of $\kappa$ becomes 
\begin{equation}
\label{eq:sme-vmf}
\hat{\kappa}_\hsme = d_n/W_n,
\end{equation}
expressed in terms of the first two sample moments of the standardized
data.  All sums here and below range over $h=1, \ldots, n$.

The score matching estimator for the von Mises
distribution can also be derived from the two trigonometric moments,
$$
E(\cos \nu \theta) = I_\nu(\kappa) / I_0(\kappa), \quad \nu \geq 0,
$$
for $\nu=1,2$.  Using the Bessel function identity
\begin{equation}
\label{eq:Bessel}
I_{\nu+1}(\kappa) = I_{\nu-1}(\kappa) - 
\frac{2 \nu}{\kappa} I_\nu(\kappa)
\end{equation}
with $\nu=1$ and simplifying yields the population version of
(\ref{eq:sme-vmf}) for the circle, $q=2$.  Analogous results using Legendre
polynomials yield the population version of (\ref{eq:sme-vmf}) for
larger values of $q$.

If $q=2$ and the data are represented in polar coordinates, $(z_{h1},
z_{h2}) = (\cos \theta_h,\ \sin \theta_h) (h=1, \ldots, n)$, then the
score matching estimators can be recast in polar coordinates.  In
particular, the estimated orientation angle for the hybrid score
matching estimator is

\begin{equation}
\label{eq:vm-mean}
\hat{\theta}_{0,\hsme} = \hat{\theta}_{0,\mle} = \atan2({\bar{S}},{\bar{C}}) =
\hat{\theta}_0, \text{  say,}
\end{equation}
where $\atan2(.)$ is defined so that $\atan2(y,x) = \theta$ if and
only if $(\cos \theta ,\sin \theta)^T \propto (x,y)^T$ for $x^2 + y^2
> 0$, and where
\begin{equation*}
\label{eq:vm-sufficient}
\bar{C}= \frac{1}{n}{\sum \cos \theta_h},\quad
\bar{S}= \frac{1}{n}{\sum \sin \theta_h}, \quad
\bar{R} = \sqrt{\bar{S} ^2 +\bar{C}^2}.
\end{equation*}
The hybrid score matching estimator of $\kappa$ can be re-expressed as
\begin{equation}
\label{eq:vm-kappa}
\hat{\kappa}_\hsme = \frac{\sum \cos(\theta_h - \hat{\theta}_0)}
{\sum \sin^2(\theta_h- \hat{\theta}_0)}  = 
\frac{n \bar{R}}{\sum \sin^2(\theta_h- \hat{\theta}_0)}.
\end{equation} 

The full score matching estimator on the circle is also
straightforward to derive, with $m=2$ and sufficient statistics $
t_1(\theta) =\cos \theta $ and $ t_2(\theta) =\sin \theta$.  Set
$\bar{C}_2 = \frac1n \sum \cos 2 \theta_h, \ \bar{S}_2 = \frac1n \sum
\sin 2 \theta_h, \ \bar{R}_2 = (\bar{C}_2^2 + \bar{S}_2^2)^{1/2}$.
Then  $W_n$ and $d_n$ have elements
$$ 
w^{(n)}_{11}= \frac{1}{2} (1-\bar{C}_2), \quad 
w^{(n)}_{12}= -\frac{1}{2} \bar{S}_2, \quad
w^{(n)}_{22}=\frac{1}{2} (1 + \bar{C}_2),
$$
and
$
d^{(n)}_1= \bar{C}, \ d^{(n)}_2 =\bar{S}$.  Since $\left|W_n\right| =
(1- \bar{R}_2^2)/4$, the full score matching estimator becomes
\begin{align}
\hat{\theta}_{0,\sme} &= \atan2\{\bar{C} \bar{S_2} + \bar{S} (1- \bar{C_2)},
\bar{C} (1+ \bar{C_2}) + \bar{S}\bar{S_2} \} \notag\\
\hat{\kappa}_\sme &=  2 \{ \bar{R}^2(1+  \bar{R_2}^2) + 
2(\bar{C}^2 - \bar{S}^2)\bar{C_2} + 4 \bar{C} \bar{S}\bar{S_2}\}^{1/2}/
(1 - \bar{R_2}^2) . \notag
\end{align} 
After rotating the data so that $\bar{S}=0$, the full score matching
estimate of $\kappa$ turns out to be the same as the hybrid estimate
if $\bar{S}_2=0$.

\end{itemize}
\subsection{Bingham distribution}
\label{sec:dda:Bingham}
\begin{itemize}
\item[(a)] The parameter matrix $A$ for the Bingham distribution in
  Table \ref{table:hybrid} is a symmetric $q \times q$ matrix, where
  without loss of generality, the trace of $A$ may be taken to be 0.
  If $A = \Gamma \Lambda \Gamma^T$ is the spectral decomposition of
  $A$, then the orthogonal matrix $\Gamma = \left[
    \gamma_{(1)},  \ldots, \gamma_{(q)} \right]$, whose columns
  are eigenvectors of $A$, represents the orientation parameters and the 
  eigenvalues $\Lambda = \diag(\lambda_j)$ represent the concentration
  parameters, with $\sum_{j=1}^q \lambda_j = \tr(A)=0$.
\item[(b)] Given the $n \times q$ original data matrix $Z$ calculate the
  moment of inertia matrix $T^{(Z)} = (1/n) Z^T Z$ and find its
  spectral decomposition $T^{(Z)} = GLG^T$.  Then the maximum
  likelihood estimate of $\Gamma$ is $G$.  Define the standardized
  data matrix by $Y = Z G$.
\item[(c)] In the reduced model, the matrix $A$ simplifies to the
  diagonal matrix $\Lambda$, with $y^T \Lambda y = \sum_{j=1}^{q-1}
  \lambda_j (y_j^2 - y_q^2)$ since $\sum_{j=1}^q \lambda_j = 0$.
\item[(d)] For the reduced model $\pi$ becomes the parameters
  $\lambda_j \ (j=1, \ldots, q-1)$ and the estimates $W_n$ and $d_n$
  have entries for $i,j = 1, \ldots, q-1$,
$$
w^{(n)}_{ij} = \begin{cases} 
\frac4n \sum \{y_{hi}^2+y_{hq}^2 -(y_{hi}^2-y_{hq}^2)^2\}, \quad i=j\\
\frac4n \sum \{y_{hq}^2 -(y_{hi}^2-y_{hq}^2)(y_{hj}^2-y_{hq}^2)\}, \quad i \neq j
\end{cases}
$$
and
$$
d^{(n)}_i = \frac{2q}{n} \sum  (y_{hi}^2 - y_{hq}^2).
$$
\end{itemize}

\subsection{Kent distribution}
\label{sec:dda:fb5}
The Fisher-Bingham distribution in Table \ref{table:hybrid} has
$q^2+3q-2$ parameters which can be estimated by the score matching
estimator.  However, this distribution has too many parameters to be
of much interest in practice.  Instead it is more useful to consider
sub-families of this distribution.  One such sub-family for $S_2$ is
the 5-parameter FB5 distribution, also known as the Kent distribution,
which forms a curved exponential family.

\begin{itemize}
\item[(a)] For the Kent distribution, $\mu = \kappa \gamma_{(1)} $ is
  assumed to be an eigenvector of $A = \Gamma \Lambda \Gamma^T$ with
  eigenvalue $\lambda_1 = 0$; the other two eigenvalues are assumed to
  be of equal size with opposite signs, $\lambda_2 = -\lambda_3 =
  \beta$.  The orthogonal matrix $\Gamma$ contains the orientation
  parameters and there are two concentration parameters, $\kappa \geq
  0$ and $\beta \geq 0$.
\item[(b)] \citet{Kent82} describes a moment estimator $\hat{\Gamma}$
  for the orientation matrix, which will also be used for the hybrid
  score matching estimator.  After standardizing to $Y = Z
  \hat{\Gamma}$, the sample mean vector becomes $\overline{y} =
  (\overline{R},0,0)^T$ and the moment of inertia matrix $T^{(Y)} =
  (1/n) Y^TY$ satisfies $(T^{(Y)})_{23} = (T^{(Y)})_{32} = 0$ and
  $(T^{(Y)})_{22} - (T^{(Y)})_{33} \geq 0$.
\item[(c)]  The reduced form of the distribution is given in Table
\ref{table:hybrid} with $\pi = (\kappa, \beta)^T$ having two components.
\item[(d)] The estimates $W_n$ and $d_n$ have entries
\begin{align*}
w^{(n)}_{11} &=  1- \sum y_{h1}^2, \quad
w^{(n)}_{12} = w_{21}^{(n)} = \frac2n \sum y_{h1}(y_{h3}^2 - y_{h2}^2), \\
w^{(n)}_{22} &= \frac4n \sum \{y_{h2}^2+y_{h3}^2 -(y_{h2}^2-y_{h3}^2)^2\}
\end{align*}
and
$$
d^{(n)}_1 = \frac2n \sum y_{h1}, \ 
d^{(n)}_2  = \frac6n \sum (y_{h2}^2 - y_{h3}^2).
$$
\end{itemize}
\subsection{Multivariate von Mises sine model on the torus}
\label{sec:dda:sine}
A general model on a product manifold $M = \prod_{r=1}^k M^{(r)} =
M^{(1)} \otimes \cdots \otimes M^{(k)}$, involving first order
interaction terms, takes the form
\begin{equation*}
\label{eq:exp-fam-prod}
f(x; \pi) \propto \exp \{ \sum_{r=1}^k \pi^{(r)T} t^{(r)}(x^{(r)}) +
\sum_{r<s} t^{(r)T}(x^{(r)}) \Omega^{(r,s)} t^{(s)}(x^{(s)}) \}
\end{equation*}
\citep{Jupp-Mardia80}.  Different sufficient statistics, of possibly
different dimensions, are allowed each manifold $M^{(r)}$.  In
principle it is also possible to include higher order interactions, 
and in some applications it may be desirable to consider submodels by
setting some of the parameters to 0.

An important example of this construction is the general multivariate
von Mises distribution on the torus $M = (S_1)^k$, for which
$t^{(r)}(z^{(r)}) = (z^{(r)}_1, z^{(r)}_2)^T = $ \\$(\cos \theta^{(r)},
\sin \theta^{(r)})^T$ comprises the two Euclidean coordinates on
each circle.  See e.g.  \citep{Mardia75,Mardia-Pat05} for the torus
case. \citet{Kume.etal13} give an extension to a product of higher
dimensional spheres, $q>2$.

Even on the torus, this model has too many interaction parameters to
be easily interpretable, so simplifications are often considered,
including a sine model and two versions of a cosine model \citep[e.g.,
][p. 501]{Mardia13}.  For this paper we limit attention to the sine
model on the torus; the cosine versions can be analyzed similarly.

\begin{itemize}
\item[(a)] The density for the sine model takes the form  in
  Table \ref{table:hybrid} where $\theta^{(r)\prime} = \theta^{(r)} -
  \theta^{(r)}_0$ \citep{Singh.etal02,Mardia.etal08}.  Here
  $\theta^{(r)}_0 \ (r=1, \ldots, k)$, denote a set of orientation or
  centering parameters.  The sine model forms a curved exponential
  family.
\item[(b)] The centering parameters can be estimated marginally by the
  sample mean directions $\hat{\theta}^{(r)}_0 \ (r=1, \ldots, k)$ on
  each circle separately.  Let $\phi_h^{(r)} = \theta_h^{(r)} -
  \hat{\theta}^{(r)}_0$ denote the standardized angles for the $n$
  data points arranged as an $n \times k$ matrix $\Phi$.  In this
  example, it is simpler to work in polar coordinates than in
  Euclidean coordinates.
\item[(c)] The reduced model takes the form
\begin{equation}
\label{eq:torus-sine-reduced}
f(\phi) = \exp\{ \sum_{r=1}^k \kappa^{(r)} \cos \phi^{(r)} +
\sum_{r<s} \lambda^{(rs)} \sin \phi^{(r)} \sin \phi^{(s)}\}
\end{equation}
and forms a canonical exponential family with $m=k + k(k-1)/2$
parameters.  The $m$-dimensional sufficient statistic, denoted
$t(\phi)$, say, can be split into two blocks, $\cos\phi^{(r)} \
(r=1, \ldots, k)$ and $\sin \phi^{(r)} \sin \phi^{(s)} \ (r<s)$.

\item[(d)] Working in polar coordinates, the $m \times k$ matrix of
  partial derivatives $\nabla^T t(\phi)$ can be similarly be
  partitioned into two blocks, where the nonzero elements are 
\begin{align*}
(\nabla \cos \phi^{(r)})_j  &=   -\sin \phi^{(r)}, \quad j=r\\
\left\{\nabla (\sin \phi^{(r)} \sin \phi^{(s)}) \right\}_j  &= 
\begin{cases} \cos \phi^{(r)} \sin \phi^{(s)}, \quad j=r,\\
              \sin \phi^{(r)} \cos \phi^{(s)},  \quad j=s,
\end{cases}
\end{align*}
for $j,r,s=1, \ldots, k, \ r<s$.   Then
$$
W_n = \frac1n \sum \nabla^T t(\phi_h) \{\nabla^T t(\phi_h)\}^T,
$$
where $\phi_h = (\phi_{h1}, \ldots, \phi_{hk})^T$ denotes the $k$-vector
of angles for row $h$ of the standardized data matrix.   

The functions $\cos \phi^{(r)}$ are eigenfunctions of the
Laplace-Beltrami operator with eigenvalue -1.  Similarly, the product
functions $\sin \phi^{(r)} \sin \phi^{(s)}$ are eigenfunctions with
eigenvalue -2.  Hence the elements of $d_n$ have entries
\begin{align*}
d^{(n)}_r &=  \frac1n \sum \cos \phi_h^{(r)} \quad (r=1, \ldots, k),\\
d^{(n)}_{rs} &= \frac2n \sum  \sin \phi_h^{(r)} \sin \phi_h^{(s)}\quad (r<s).
\end{align*}
\end{itemize}

\begin{figure}[h!]  
\centering 
\includegraphics[width=6in, height=4in]{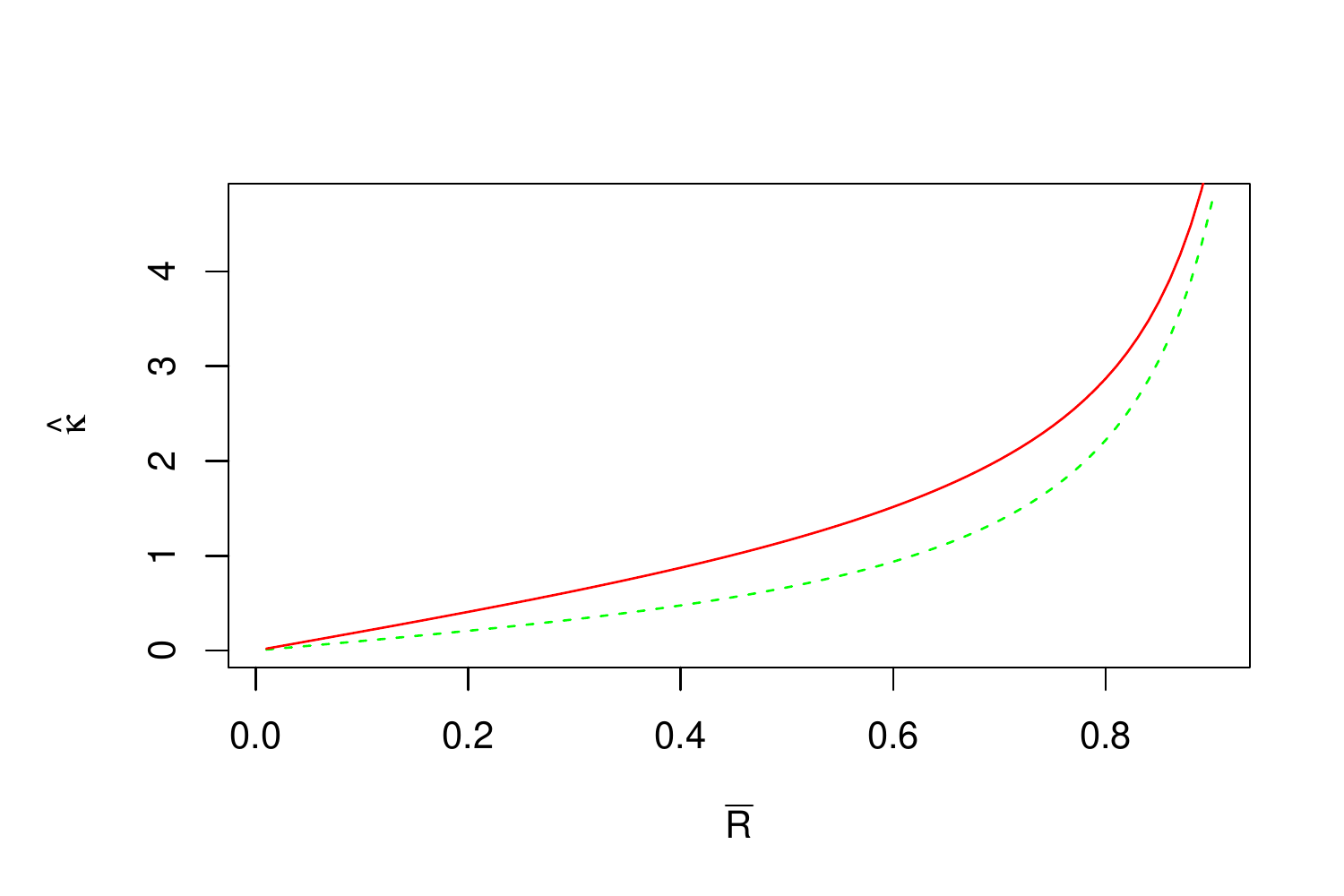}
\caption{Comparison of two estimators of $\kappa$ for the von Mises
  distribution with $n=2$: maximum likelihood estimator (solid line)
  and score matching estimator (dashed line).}
\label{fig:vm-kappa-n2}
\end{figure}

\begin{figure}[h!]  
\centering 
\includegraphics[width=6in, height=4in]{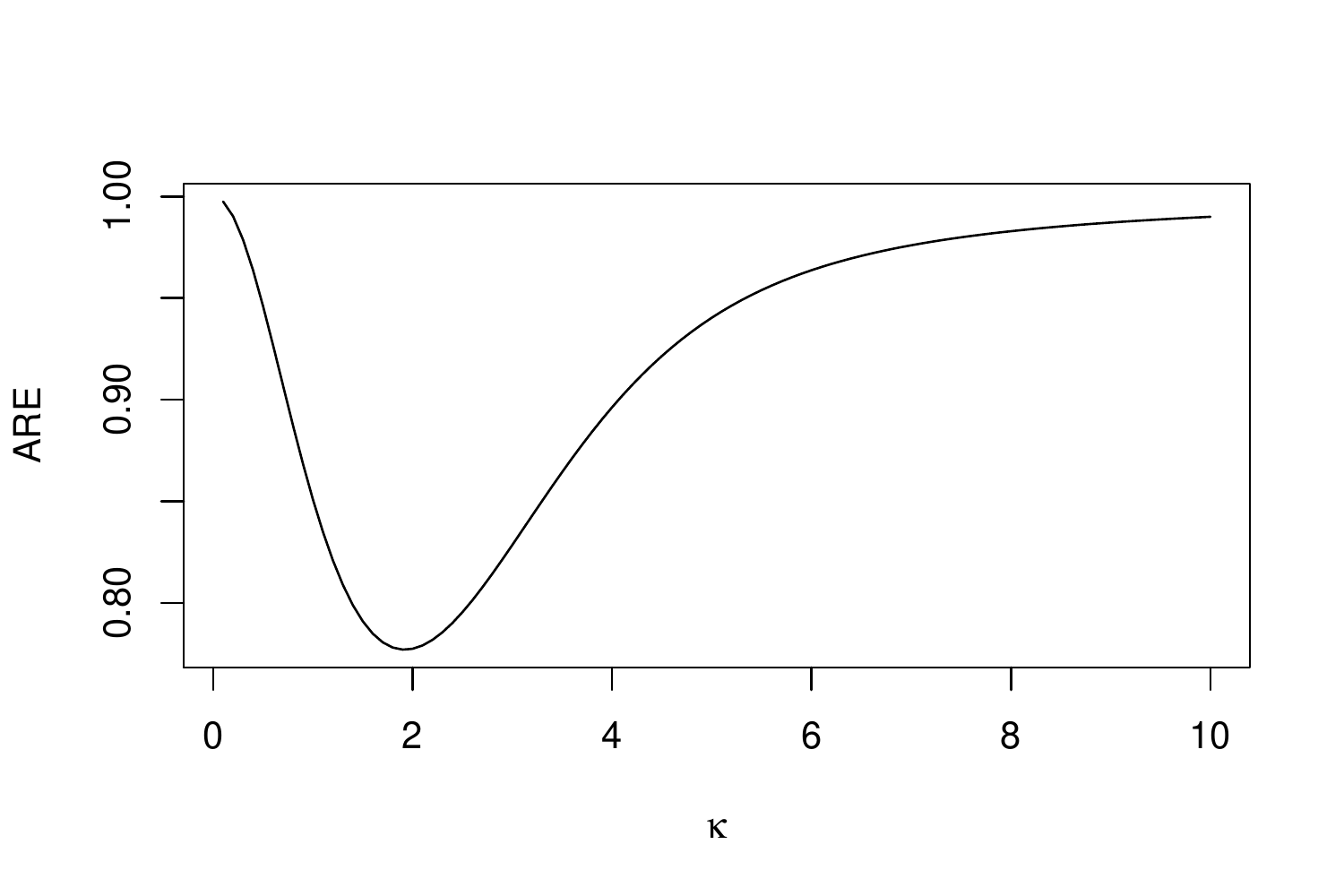}
\caption{Asymptotic relative efficiency (ARE), as a ratio of mean
  squared errors, for the estimation of $\kappa$ for the von Mises
  distribution, comparing the score matching estimator to the maximum
  likelihood estimator.}
\label{fig:are}
\end{figure}

\section{Efficiency study for the von Mises distribution}
\label{sec:efficiency}

\begin{table}
  \caption{Relative efficiency, as a ratio of mean squared errors,
    comparing the score matching estimator of $\log \kappa$ to the maximum
    likelihood estimator.  The final column ARE gives the asymptotic 
    relative efficiency as $n \rightarrow \infty$, taken from (\ref{eq:are}).}

\label{table:eff}
\begin{center}
\begin{tabular}{rrrrrr}	
$\kappa$ & n=2 &   $n=10$  &  $n=20$   & $n=100$ &  ARE\\                                                         
0.5 &	100 &	89 &	93 &	95 & 95 \\
1   &	 98 &	88 &	86 &	85 & 85 \\ 
2   &    98 &   92 &	85 &	79 & 78 \\  
10  &   100 &  100 &	99 &	99 & 99 \\ 
\end{tabular}\end{center}
\end{table}
For the von Mises distribution on the circle, it is possible to study
the behaviour of the score matching estimator in more detail, both
empirically and analytically.  The score matching estimators were
described in (\ref{eq:vm-mean}) and (\ref{eq:vm-kappa}).  The maximum
likelihood estimator of $\theta_0$ is the same as hybrid score matching
estimator.  The maximum likelihood estimator of $\kappa$ is
$\hat{\kappa}_{MLE} = A_1^{-1}(\bar{R}),$
where $A_\nu(\kappa) = I_\nu(\kappa)/I_0(\kappa), \ \nu \geq 0$.

Table \ref{table:eff} gives the relative efficiency, as a ratio of
mean squared errors, comparing the score matching estimator of $\log
\kappa$ to the maximum likelihood estimator.  Each entry is based on
100,000 simulated datasets.  The log transformation is used to improve
the numerical stability of the results, though in the final column for
asymptotic relative efficiency, the use of a transformation makes no
difference.  The simulation from von Mises distribution was carried
out in R \citep{R} using the package CircStat
\citep{CircStats12,Jamm-Sen01}.  For all cases, the relative
efficiency is at least 78\% and is generally much closer to 100\%.

The case $n=2$ is interesting because it is possible to write the
estimators in closed form and so compare their behaviour in detail.  Further,
the hybrid and full score matching estimators are identical in this
setting.  After centering, the standardized data take the form $\pm
\theta$ for a single value of $\theta, \ 0<\theta<\pi/2$, where for
simplicity we exclude the extreme possibilities $\theta=0,\pi/2$.
Then
\begin{equation}
\label{eq:vm:n=2}
\hat{\kappa}_\sme = \cos\theta / \sin^2 \theta, \quad
\hat{\kappa}_\mle = A_1^{-1}(\cos \theta).
\end{equation}
The two estimators are compared in Figure \ref{fig:vm-kappa-n2} as a
plot of the estimated $\kappa$ vs. $\overline{R} = \acos \theta$.
Although (\ref{eq:vm:n=2}) does not provide enough information to
compute the relative efficiency, even in the extreme
setting $n=2$ the two estimators are reasonably similar, with the
difference tending to 0 as $\overline{R} \rightarrow 0$ and with the
relative difference tending to 0 as $\overline{R} \rightarrow 1$. The
maximum difference between the two estimators is about
$\hat{\kappa}_\mle - \hat{\kappa}_\sme = 2.46-1.80 = 0.66$ when
$\overline{R}=0.76$, i.e. $\theta = 41^o$.

In the limiting case $n \rightarrow \infty$, it is possible to compute
analytically the asymptotic relative efficiency of the hybrid score
matching estimator for $\kappa$, relative to the maximum likelihood
estimator, 
\begin{equation}
\label{eq:are}
ARE = \frac {A_1^2(\kappa)} { \left\{2\kappa - 3 A_1(\kappa) \right\}
\left\{ \kappa - \kappa A_1^2(\kappa) - A_1(\kappa) \right\} }. 
\end{equation}
A plot of (\ref{eq:are}) is given in Figure \ref{fig:are} and a proof is 
given in the Appendix.

\section{Numerical examples}
\label{sec:examples}
\subsection{Kent distribution}
\label{sec:examples:fb5}
\begin{table}
  \caption{Great Whin Sill data on the sphere $S_2$, with sample size 
    $n=34$, and fitted by the Kent distribution.   Estimates of $\kappa$ 
    and $\beta$ are given by  four different methods described in the text.}
\label{table:fb5}
\begin{center} \begin{tabular}{ccc}
    Method & $\hat{\kappa}$ & $\hat{\beta}$ \\
    Hybrid  SME      &  42.13 & 9.34\\
    Hybrid MLE &  42.16 &  9.27\\
    Hybrid approximate MLE &  41.76 &  8.37\\
    Full MLE   &  42.41 &  9.28
\end{tabular} \end{center} \end{table}

The Great Whin Sill dataset was analyzed in \citet{Kent82} and is
presented here to illustrate various estimates of the concentration
parameters for the Kent distribution; see Table \ref{table:fb5}.  The
first three methods of estimation are hybrid estimators.  Hence, as
discussed in Section \ref{sec:dda:fb5}, the moment estimator is used
for the orthogonal matrix representing the orientation
parameters. After rotation of the data, the reduced model involves
just the two concentration parameters.

Here are further details about the estimators in Table
\ref{table:fb5}.  The hybrid score matching estimator for $\kappa$ and
$\beta$ was described in Section \ref{sec:dda:fb5}.  The ``hybrid
maximum likelihood estimate'' is the maximum likelihood estimate for
the concentration parameters after rotating the data using the moment
estimate of orientation (eqn. (4.8) in \citet{Kent82}). The ``hybrid
approximate maximum likelihood estimate'' is the same, but using a
normal approximation for the normalizing constant (eqn. (4.9) in
\citet{Kent82}).  The ``full maximum likelihood estimate'' involves
maximizing the 5-parameter likelihood over both the orientation and
concentration parameters; the estimated orientation is negligibly
different from the moment estimator and is not reported here.  All the
estimates of the concentration parameters are close together.

\subsection{Torus}
\label{sec:examples:sine}

\begin{table}
  \caption{Isoleucine data, component 1, on the torus $(S_1)^4$, with
    sample size $n=23$, and fitted using the multivariate von Mises sine
    model.  Estimates of the concentration parameters are given by three
    different methods: the hybrid score matching estimator (SME), the 
    hybrid composite likelihood estimator and the hybrid approximate
    maximum likelihood estimator (AMLE).}
\label{table:torus}

\begin{center}
\begin{tabular}{rrrrcrrrrcrrrr}
\multicolumn{4}{c}{SME} &$\ \ $& \multicolumn{4}{c}{composite} & $\ \ $&
\multicolumn{4}{c}{AMLE}  \\
4.3 & -6.3 & 3.4 & 2.5  &&  4.9 &-7.3 & 4.0&  3.4  && 6.6 & -5.3 & 3.1 & 3.1 \\
*   & 45.4 &  24.5 & 6.03 && * & 48.4& 26.3 &  7.1 && * & 47.1 & 24.4 & 6.0\\
*    & *   &  59.3 & -4.6 && * & * & 61.6 & -5.6 && * & * &      60.1 & -4.3\\
* & * & * &           2.2 && * & * & * &     2.5 && * & * & * & 3.6 \\
\end{tabular}
\end{center}
\end{table}

In \citet{Mardia.etal12}, $k=4$ angles from the amino acid isoleucine
were modelled by a mixture of multivariate von Mises sine models and
grouped into 17 clusters.  Here we look at just one of those clusters, Cluster 1, 
and look at the fits to the concentration parameters from three
estimation methods. All the methods are hybrid methods.  Thus in each
case location is estimated using the moment estimator given by
the sample mean directions for each of the $k=4$ angles.  After
rotation the reduced model involves just the $k(k+1)/2 = 10$
concentration parameters.

The hybrid score matching estimator was described in Section
\ref{sec:dda:sine}.  The hybrid approximate maximum likelihood estimator uses
a high concentration normal approximation for the normalizing constant
in the reduced model; the numerical values were given in
\citet{Mardia.etal12}.  The hybrid composite likelihood estimator is based on
the conditional von Mises distribution in the reduced model for each
angle given the remaining angles; the methodology is summarized in
\citet{Mardia.etal09}.

The estimated parameters are given in Table \ref{table:torus}.  The
estimates are presented as a symmetric matrix: the diagonal elements
are the $\kappa^{(r)}$ and the off-diagonal elements are the
$\lambda^{(rs)}$ in (\ref{eq:torus-sine-reduced}), given for clarity
just in the upper triangle.  In general all the estimates match
reasonably closely.

\section{Discussion}
\label{sec:future}

Methods of estimation for exponential families on manifolds can be
divided into at least three broad categories:
\begin{itemize}
\item Maximum likelihood estimators, both the exact version and
  approximate versions.  Although the exact version is preferred in
  principle, there may be problems in practice evaluating the
  normalizing constant.  Hence approximations may be used, such as (a)
  saddlepoint \citep{Kume.etal13}, (b) holonomic \citep{Sei-Kume15},
  and (c) approximate normality under high concentration.

\item Composite maximum likelihood estimation. Suppose a point on the
  manifold can represented as a set of variables, such that the
  conditional distribution of each variable given the rest is
  tractable. Then the composite likelihood is the product of
  the conditional densities.  In some cases this method can be very
  efficient \citep{Mardia.etal09}.
 
\item Score matching estimators.  As shown in this paper these
  estimators often reduce to the solution to a set of linear equations
  based on sample moments of the data.  Hence the method is easy to
  implement and straightforward to apply to large datasets, including
  streaming data.
\end{itemize}
Strictly speaking, an approximate maximum likelihood estimator, at a
fixed level of approximation, will not be consistent as the sample
size $n \rightarrow \infty$. The score matching estimator is always
consistent under the mild regularity conditions (A1)--(A2) on $f$ in
Section \ref{sec:sme-derivation}, with an asymptotic variance at
least as large as the maximum likelihood estimator.  The the numerical
examples here suggest the efficiency of the hybrid score
matching estimator, compared to the maximum likelihood estimator, will
often be close to 1.

Score matching estimators can also be developed for distributions on
noncompact manifolds, including Euclidean spaces.  The simplest
example is  the multivariate normal distribution \citep{Hyv05} where it turns
out that the score matching estimator is identical to the maximum
likelihood estimator.  Many directional distributions are
approximately normal under high concentration.  Hence we expect high
efficiency of the score matching estimator in this setting.  See
Table \ref{table:eff} for confirmation in the von Mises case.

In the directional setting, the score matching estimator can often be
interpreted as a ``double moment estimator''.  For example, for the
von Mises-Fisher distribution the sufficient statistic is a linear
function of $z$, but the matrix $W_n$ in the score matching estimator
involves quadratic functions of $z$.  Similarly, for the Bingham
distribution the sufficient statistic is a quadratic function of $z$,
but the matrix $W_n$ in the score matching estimator involves quartic
functions of $z$.

This paper has emphasized the setting where $M$ is a sphere or a
product of spheres.  Work is in progress to investigate the score
matching estimator for models on other important manifolds such as
Stiefel and Grassmann manifolds.

\section*{Acknowledgements}
The authors are grateful to Derek Harland, John Wood and Peter Kim for
helpful discussions about the differential geometry, and to Peter
Forbes for helpful discussions of the von Mises case.  The first author
also thanks the Indian Institute of Management Ahmedabad for 
hospitality during his visits to carry out some of this work.

\section{Appendix}
This section gives a proof of the asymptotic relative efficiency for
the score matching estimator of $\kappa$ in (\ref{eq:are}).  For
simplicity the proof focuses on the estimator in the reduced model.

First, it can be shown that, asymptotically,
\begin{equation}
\label{eq:sme-asymp}
n \var(\hat{\kappa}_\hsme) = 
\frac{\kappa}{A_1^2(\kappa)} \{2\kappa - 3 A_1(\kappa)\}.  \tag{A1} 
\end{equation} 
To verify this equation note that, treating $\theta$ as random from
the von Mises distribution, $\var(\cos \theta) = 1 - A_1^2(\kappa) -
A_1(\kappa)/\kappa,$ $\var(\sin^2 \theta) = \frac{1}{8} \left\{ 3 +
  A_4(\kappa) - 4 A_2(\kappa) \right\} - A_1^2(\kappa)/\kappa^2$ and
$\cov(\cos \theta, \sin^2 \theta) = \frac{1}{4} \left\{ A_1(\kappa) -
  A_3(\kappa) \right\} - A_1^2(\kappa) / \kappa$.  Using the delta
rule gives
$$ 
\var(\hat{\kappa}_\hsme) = \frac{\kappa}{8nI_{1}A_1(\kappa)}
\left\{\kappa (8 + 3 \kappa^2) I_{0} - 8 I_{1} - 4
  \kappa^2 I_{1} - 4 \kappa^3 I_2 + 4 \kappa^2
  I_{3} + \kappa^3 I_{4}\right\},
$$ 
with the shorthand notation $I_\alpha = I_\alpha(\kappa)$.  Repeated
use of (\ref{eq:Bessel}) leads to (\ref{eq:sme-asymp}).

For the maximum likelihood estimator it can be shown that,
asymptotically,
\begin{equation}
\label{eq:mle-asymp}
 n \ \var( \hat{\kappa}_\mle) = \frac{1}{1 - A_1^2(\kappa) - 
A_1(\kappa)/\kappa}. \tag{A2}
\end{equation} 
For both estimators the asymptotic variance converges to 2 as $\kappa
\rightarrow 0$ and is asymptotic to $2 \kappa^2$ as $\kappa
\rightarrow \infty$.  Hence for small and large $\kappa$ , the
asymptotic relative efficiency tends to 1.  Combining
(\ref{eq:sme-asymp}) and (\ref{eq:mle-asymp}) yields (\ref{eq:are}).

The discussion here has emphasized the hybrid score matching estimator
and the maximum likelihood estimator for the full model with two
unknown parameters.  But in fact there are other models and estimators
to consider, including three versions of the score matching
estimator and two versions of the maximum likelihood estimator.  The
hybrid score matching estimator and the full score matching estimator
are defined for the full model with two unknown parameters $\theta_0$
and $\kappa$. The score matching estimator can also be defined under
the reduced model with just one unknown parameter.  It turns out that
all three estimators of $\kappa$ have the same asymptotic variance.
Similarly, the maximum likelihood estimator of $\kappa$ can be defined
in the setting of the full model or the reduced model.  Again both
estimators have the same asymptotic variance.  Hence the asymptotic
relative efficiency takes the same value for all these possibilities.

\bibliographystyle{apalike}
\bibliography{refs-sme} \
\end{document}